\documentclass[twoside, 12pt]{report}
\usepackage[english]{babel}
\usepackage{epsfig}
\usepackage{psfig}
\usepackage{amssymb,amsmath,amsfonts,soul,amsthm,enumerate}
\usepackage{graphicx}
\usepackage[cp1251]{inputenc}
\usepackage[T2A]{fontenc}
\usepackage{mathrsfs}
\usepackage{cite} 

\usepackage{tikz}
\newcommand{\tikzmark}[1]{\tikz[overlay,remember picture] \node (#1) {};}
\tikzset{square arrow/.style={
		to path={-- ++(0,-.45)  -| (\tikztotarget) \tikztonodes},below,pos=.25}}

\def\@evenfoot{}
\def\@oddfoot{}

\usepackage{authblk}

\begin{document}

\def\@evenhead{\vbox{\hbox to \textwidth{\thepage\leftmark}\strut\newline\hrule}}

\def\@oddhead{\raisebox{0pt}[\headheight][0pt]{%
\vbox{\hbox to \textwidth{\rightmark\thepage\strut}\hrule}}}

\def\bibname{\vspace*{-30mm}{\centerline{\normalsize References}}}

\newpage
\normalsize
\thispagestyle{empty}

\vskip 5 mm

\centerline{\bf NUMERICAL COMPARISON OF ITERATIVE AND}
\centerline{\bf FUNCTIONAL-ANALYTICAL ALGORITHMS FOR}
\centerline{\bf INVERSE ACOUSTIC SCATTERING}

\vskip 0.3cm
\centerline{\bf A.S. Shurup$^{ \textbf{1),2)}}$}

\vskip 0.3cm

\begin{center}
\begin{minipage}{118mm}
	\small
	\textit{$^{1)}$~M.V.~Lomonosov Moscow State University, Faculty of Physics, Acoustics Department,
		Leninskie Gory, Moscow 119991, Russia}.
		
	\textit{$^{2)}$~Sсhmidt Institute of Physics of the Earth
		of the Russian Academy of Sciences, 
		B. Gruzinskaya str., 10, build. 1, Moscow 123242, Russia.}

	\centerline{ \textit{e-mail: shurup@physics.msu.ru}}

\end{minipage}
\end{center}

\vskip 0.3cm
\noindent{\bf Abstract}
{\small In this work the numerical solution of acoustic tomography problem based on the iterative and functional-analytical algorithms is considered. The mathematical properties of these algorithms were previously described in works of R.G.~Novikov for the case of the Schr{\"o}dinger equation. In the present work, for the case of two-dimensional scalar Helmholtz equation, the efficiency of the iterative algorithm in reconstruction of middle strength scatterers and advantages of the functional-analytical approach in recovering strong scatterers are demonstrated. A filtering procedure is considered in the space of wave vectors, which additionally increases the convergence of the iterative algorithm. Reconstruction results of sound speed perturbations demonstrate the comparable noise immunity and resolution of the considered algorithms when reconstructing middle strength scatterers. A comparative numerical investigation of the iterative and functional-analytical algorithms in inverse acoustic scattering problems is implemented in this work for the first time.
}
\vskip 0.3cm

\setcounter{figure}{0}

\renewcommand{\thesection}{\large 1}

\section{\large Introduction}

Acoustic waves have a unique penetrating ability and can propagate in almost any natural environments. This caused the rapid development of various "remote" or "non-invasive" methods of acoustic diagnostics of natural mediums, among which acoustic tomography plays a special role. The most commonly used areas of acoustic tomography applications include medical surveys \cite{label1,label2}, tomography of an inhomogeneous moving ocean \cite{label3,label4}, as well as seismoacoustic methods of investigating the Earth \cite{label5,label6}. Usually, the solution of tomographic problems in these practical applications is based on a linear approximation, which leads to simple relationships between the properties of medium and parameters of acoustic signals that have passed through this medium. In the case of a strong deviation from linear approximation, multistep iterative procedures and various regularization methods, based on a priori information, are used to refine reconstruction results. From a mathematical point of view, the acoustic tomography is a special case of a more general class of inverse scattering problems. There are known mathematically rigorous functional-analytical methods \cite{label7,label8,label9,label9_1,label10,label11,label12,label13} for solving inverse scattering problems in quantum mechanical applications without iterations and additional regularization procedures. Numerical modeling of these algorithms has shown its perspectives for solving acoustic tomography problems \cite{label14,label15}. At the same time, the numerical study of the multichannel version of the functional algorithm \cite{label12} for the purposes of ocean mode tomography \cite{label16} has revealed that for experimental conditions it is not always possible to obtain reconstruction results with the desired accuracy. This required a further search for methods of inverse problem solution in this case. The solution of inverse dynamic problems taking into account the interaction of propagating waves, which is an analogue of multichannel scattering, was considered in works devoted to the development of boundary control methods \cite{label17,label18}. A perspective solution for multichannel inverse problems was proposed by R.G. Novikov in \cite{label19}, where an iterative algorithm for solving the inverse scattering problem for the Schr{\"o}dinger equation at a fixed energy was considered. In contrast to rigorous functional-analytical methods \cite{label7,label8,label9,label9_1,label10,label11,label12,label13}, this approach gives an approximate solution, but it is very flexible and can be adapted to solve various problems of acoustic tomography. The known results of numerical simulation of this approach in quantum mechanical applications \cite{label20} demonstrate its fast convergence and acceptable noise immunity. These results also indicate the perspectives of this algorithm in problems of acoustic tomography. It should be noted that in \cite{label19} mathematically rigorous estimates of convergence for the considered iterative algorithm were obtained that distinguishes it from other known to author iterative algorithms for solving inverse problems of acoustic scattering \cite{label2}, convergence conditions of which were studied mainly qualitatively, based on physical considerations.

In the present work, capabilities of the iterative algorithm \cite{label19} in problems of acoustic tomography are analyzed. The numerical simulation of the sound speed reconstruction is considered. The reconstruction results are compared with the estimates obtained by the functional-analytical algorithm \cite{label21,label22}, capabilities of which in the problems of acoustic tomography have already been studied before \cite{label23,label24}.

\renewcommand{\thesection}{\large 2}
\section{\large Statement of inverse problem. Recalculation of fields measured at a boundary of tomography domain to scattering amplitudes}

The two-dimensional inverse problem is considered. The solution of two-dimensional problems is required, for example, in ocean adiabatic mode tomography \cite{label4,label25}, when the three-dimensional problem of reconstructing ocean parameters is approximated by a set of independent two-dimensional inverse problems for individual hydroacoustic modes. The two-dimensional tomographic problems are also common in medical applications \cite{label26,label27,label28}. It should be noted that the initial relations of the considered algorithm \cite{label19} are given for the space dimension greater or equal to 2, which allows, if necessary, to write out the relations for solving the three-dimensional problem of acoustic tomography.

\begin{figure}[t!]
	\centerline{\epsfig{file=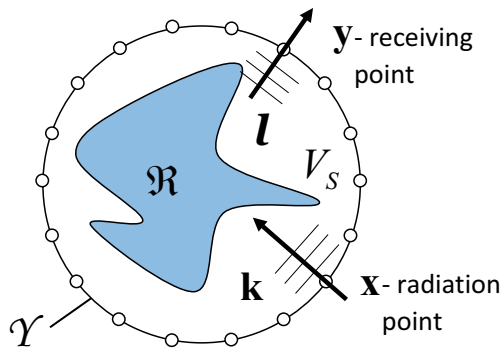, width = 7 cm}}
	\caption{Tomography area $V_S$ contains the scattering region  $\mathfrak{R}$; locations of quasi-point sources and receivers at the boundary $\Upsilon$ are described by vectors $\mathbf{x}$ and $\mathbf{y}$, respectively; wave vectors $\mathbf{k}$ and $\boldsymbol{\ell}$ show the directions of incident and scattered plane waves.}
	\label{figure_1}
\end{figure}

It is assumed that on the boundary $ S $  of tomography area $ V_S $  there are transducers, which are equivalent to point ones (such transducers will be referred below as quasi-point transducers), emitting and receiving acoustic fields. In the general case, the emitters are located on some boundary $ \mathcal{X}$ at points with radius vectors $ \mathbf{x}\in\mathcal{X} $  (where  $ \mathcal{X} \cap V_S = \varnothing $), receivers are located on a boundary  $\Upsilon $ at points with radius vectors $ \mathbf{y}\in\Upsilon $  (where $ \Upsilon \cap V_S = \varnothing $). 
However, below it is assumed that the boundary of emission $ \mathcal{X}$ and the boundary of reception $\Upsilon $ coincide with each other and with the border $ S $ of region $ V_S $: $\mathbf{x}, \mathbf{y} \in \Upsilon $, where $ \mathcal{X} \equiv \Upsilon \equiv S$, although this condition is not principal. 
Inside the region $ V_S $ there is an inhomogeneity (scatterer), which is nonzero only inside the scattering domain $\mathfrak{R}$, which lies entirely inside $ V_S $: $\mathfrak{R} \subset V_S$ (figure~\ref{figure_1}). 
For simplicity, an inhomogeneity will be considered, which is described by a sound speed perturbation only. A more general case of joint reconstruction of sound speed, density, absorption coefficient and flows can also be implemented on the basis of relations given below together with results obtained earlier in \cite{label29}.

The spatial distribution of complex spectral amplitude of acoustic pressure $ p(\mathbf{r}; \omega_j) $ in the considered inhomogeneous medium, which is characterized by the sound speed  $c(\mathbf{r})$ and a constant value of density with negligibly small influence of absorption and flows, is described by the equation \cite{label30}:
\begin{equation}
\label{Helm_eq}
\nabla^2 p(\mathbf{r}; \omega_j) + k^2_{0j} \ p(\mathbf{r}; \omega_j) = v(\mathbf{r}, \omega_j) \  p(\mathbf{r}; \omega_j) ,
\end{equation} 
\noindent
where in the considered case 
\begin{equation}
\label{Scatt_func}
v(\mathbf{r}, \omega_j) = \omega^2_j \left(\frac{1}{c^2_0} - \frac{1}{c^2(\mathbf{r})}\right) ,
\end{equation} 

\noindent
here $\mathbf{r}$ is a radius vector, $\omega_j$ is a circular frequency, $c_0$ and $k_{0j} = \omega_j \big/ c_0$ are the sound speed and the wavenumber in a background, stationary, non-absorbing medium; the time dependence is assumed in the form $\sim \text{exp}(- i \omega_j t)$. 
The index $j$ was introduced to note different frequencies $\omega_j$, which are used further to describe the polychromatic sounding regime when working with noisy scattering data. 
Because the sources and receivers are assumed to be quasi-point transducers, the acoustic fields at an arbitrary point $\mathbf{r}$ are the Green's functions  $G_0(\mathbf{r}, \mathbf{x}; \omega_j)$ or $G(\mathbf{r}, \mathbf{x}; \omega_j)$ in the absence or in the presence of a scatterer, respectively.
In the considered two-dimensional case $G_0(\mathbf{r}, \mathbf{x}; \omega_j) = - (i/4) \, H^{(1)}_0(k_{0j} |\mathbf{r} - \mathbf{x}|)$, where $ H^{(1)}_0 $ is the Hankel function of the zero order and the first kind. 
It is required to reconstruct the scatterer $v(\mathbf{r}, \omega_j), \mathbf{r} \in V_S$, by using the experimentally measured fields, i.e. scattering data $G(\mathbf{r}, \mathbf{x}; \omega_j)$. 
The described statement of problem is quite general for two-dimensional inverse problems and was considered earlier in \cite{label14, label15} for modeling various versions of the functional-analytical algorithms \cite{label12, label21, label22}.

It should be noted that in the considered iterative algorithm \cite{label19}, the initial data are not the total fields, but the scattering amplitudes, which characterize the scattered fields in the far zone \cite{label14,label15}.
When a plane wave $u_0(\mathbf{z}, \mathbf{k}; \omega_j) = \text{exp}(i \mathbf{k} \mathbf{z})$ (here $\mathbf{k}$ is the wave vector of incident field, $\mathbf{z}$ is an arbitrary point) is scattered by the inhomogeneity, the total field $u(\mathbf{z}, \mathbf{k}; \omega_j)$ in the far zone has the asymptotic \cite{label14}: 
\begin{equation}
\label{field_far_zone}
\begin{split}
& u(\mathbf{z}, \mathbf{k}; \omega_j) = \\
& = \text{exp}(i \mathbf{k} \mathbf{z}) + C_D \frac{\text{exp} \left(i k_{0j} | \mathbf{z} |\right)} {\sqrt{| \mathbf{z} |^{D-1}}} \ f(\mathbf{k},  \boldsymbol{\ell} = k_{0j} \frac{\mathbf{z}}{|\mathbf{z}|}; \omega_j) + o \left( \frac{1}{\sqrt{| \mathbf{z} |^{D-1}}}\right), \ | \mathbf{z} | \rightarrow \infty,
\end{split}
\end{equation} 


\noindent
where $C_D = -\pi \sqrt{\pi} (1+i) \big/ \sqrt{k_{0j}} $ if $D = 2$ ; $C_D = - 2 \pi^2 $ if $D = 3$; $D$ is the space dimension. 
The relation (\ref{field_far_zone}) includes the scattering amplitude $f(\mathbf{k}, \boldsymbol{\ell}; \omega_j)$, where $\boldsymbol{\ell}$ is the wave vector of scattered wave in the far zone; $\mathbf{k}, \boldsymbol{\ell} \in \mathbb{R}^2$, $\mathbf{k}^2 = \boldsymbol{\ell}^2 = k^2_{0j}$. 
As a result, the numerical implementation of the iterative algorithm \cite{label19} requires recalculation of scattering data in the form of acoustic fields $G_0(\mathbf{r}, \mathbf{x}; \omega_j)$, $G(\mathbf{r}, \mathbf{x}; \omega_j)$, emitted and received by quasi-point transducers, into scattering amplitudes $f(\mathbf{k}, \boldsymbol{\ell}; \omega_j)$, i.e. as if the tomography area was sounded by plane waves and the scattered fields were recorded by plane transducers. 
Such a recalculation of data measured in the near field into data as if they were measured in the far field is of independent particular interest, since many algorithms for solving inverse problems \cite{label7,label8,label9,label9_1,label10,label11,label12,label13} use the scattering amplitudes rather than total fields as input data. 
A possible way of such a recalculation is described in \cite{label14}. 
Below, the simplified version of calculation $f(\mathbf{k}, \boldsymbol{\ell}; \omega_j)$ based on $G_0(\mathbf{r}, \mathbf{x}; \omega_j)$, $G(\mathbf{r}, \mathbf{x}; \omega_j)$ is presented, which requires fewer numerical operations in comparison with \cite{label14} that is important for improving the performance of entire iterative algorithm.

The stages of initial algorithm \cite{label14} for recalculating fields $G_0(\mathbf{r}, \mathbf{x}; \omega_j)$, $G(\mathbf{r}, \mathbf{x}; \omega_j)$ into scattering amplitudes $f(\mathbf{k}, \boldsymbol{\ell}; \omega_j)$, which are valid for both two-dimensional $D = 2$ and three-dimensional $D = 3$ problems, are as follows \cite{label11, label14}:
\begin{enumerate}[1)]
	
	\item Finding the function $(\mathscr{F} - \mathscr{F}_0)(\mathbf{y^\prime},\mathbf{y^{\prime\prime}}; \omega_j)$, which characterizes the kernel of the Dirichlet-Neumann operator \cite{label11}, from a system of linear integral equations written for different positions of sources and receivers:	
	\begin{equation}
	\label{F_F0_int}
		\begin{split}
			& \int\limits_\Upsilon d \mathbf{y^\prime} \int\limits_\Upsilon d \mathbf{y^{\prime\prime}} \ G_0(\mathbf{y} - \mathbf{y^\prime}; \omega_j) (\mathscr{F} - \mathscr{F}_0)(\mathbf{y^\prime},\mathbf{y^{\prime\prime}}; \omega_j) G(\mathbf{y^{\prime\prime}}, \mathbf{x}; \omega_j) = \\
			& = G(\mathbf{y}, \mathbf{x}; \omega_j) - G_0(\mathbf{y} - \mathbf{x}; \omega_j), \ \mathbf{x}, \mathbf{y} \in \Upsilon.
		\end{split}
	\end{equation} 
	
	\item Calculation of intermediate function $A(\mathbf{y},\mathbf{y^{\prime\prime}}; \omega_j)$, $\mathbf{y^{\prime\prime}} \in \Upsilon$:
	\begin{equation}
	\label{A_int}
		A(\mathbf{y},\mathbf{y^{\prime\prime}}; \omega_j) = \int\limits_\Upsilon \ G_0(\mathbf{y} - \mathbf{y^\prime}; \omega_j) (\mathscr{F} - \mathscr{F}_0)(\mathbf{y^\prime},\mathbf{y^{\prime\prime}}; \omega_j) d \mathbf{y^\prime}.
	\end{equation} 
	
	\item Solution of the system of integral equations relative to the field $u(\mathbf{y}, \mathbf{k}; \omega_j)$ occurring at the reception points $\mathbf{y} \in \Upsilon$ as a response to the incident plane wave $u_0(\mathbf{y}, \mathbf{k}; \omega_j) = \text{exp}(i \mathbf{k} \mathbf{y})$:
	\begin{equation}
	\label{u_int}
		u(\mathbf{y}, \mathbf{k}; \omega_j) = \text{exp}(i \mathbf{k} \mathbf{y}) + \int\limits_\Upsilon \ A(\mathbf{y}, \mathbf{y^{\prime\prime}}; \omega_j) u(\mathbf{y^{\prime\prime}},\mathbf{k}; \omega_j) d \mathbf{y^{\prime\prime}}, \ \mathbf{y} \in \Upsilon.
	\end{equation} 
	
	\item Calculation of scattering amplitude from the functions $(\mathscr{F} - \mathscr{F}_0)(\mathbf{y^\prime},\mathbf{y^{\prime\prime}}; \omega_j)$, $u(\mathbf{y}, \mathbf{k}; \omega_j)$:
	\begin{equation}
	\label{f_int}
	\begin{split}
		f(\mathbf{k}, \boldsymbol{\ell}; \omega_j) = \frac{1}{(2 \pi)^D} \int\limits_\Upsilon d \mathbf{y^\prime} \int\limits_\Upsilon d \mathbf{y^{\prime\prime}} \ \text{exp}(- i \boldsymbol{\ell} \mathbf{y^\prime}) & (\mathscr{F} - \mathscr{F}_0)(\mathbf{y^\prime},\mathbf{y^{\prime\prime}}; \omega_j) u(\mathbf{y^{\prime\prime}},\mathbf{k}; \omega_j), \\
		& D = 2, 3.
	\end{split}
	\end{equation} 	
	
\end{enumerate}


The calculation of scattering amplitude $f(\mathbf{k}, \boldsymbol{\ell}; \omega_j)$ can be implemented not in four steps (\ref{F_F0_int})-(\ref{f_int}) described above, but in only two steps:
\begin{enumerate}[1)]
	\item Finding the intermediate function $\Phi(\mathbf{y^\prime},\mathbf{y^{\prime\prime}}; \omega_j)$, from the system of integral equations:
		\begin{equation}
		\label{Phi_int}
			\begin{split}
			& \int\limits_\Upsilon d \mathbf{y^\prime} \int\limits_\Upsilon d \mathbf{y^{\prime\prime}} \ G_0(\mathbf{y} - \mathbf{y^\prime}; \omega_j) \Phi(\mathbf{y^\prime},\mathbf{y^{\prime\prime}}; \omega_j) G_0(\mathbf{y^{\prime\prime}}, \mathbf{x}; \omega_j) = \\
			& = G(\mathbf{y}, \mathbf{x}; \omega_j) - G_0(\mathbf{y} - \mathbf{x}; \omega_j), \ \mathbf{x}, \mathbf{y} \in \Upsilon,
		\end{split}
		\end{equation} 
		\noindent
		where in contrast to (\ref{F_F0_int}), for the integration only the Green's functions $G_0$ are used, which can be calculated analytically.
		
		\item Calculation of scattering amplitude using the double Fourier transform (taking into account the signs in the exponential for different arguments) of the function $\Phi(\mathbf{y^\prime},\mathbf{y^{\prime\prime}}; \omega_j)$ on the boundary $\Upsilon$:
		\begin{equation}
		\label{f_int_new}
		f(\mathbf{k}, \boldsymbol{\ell}; \omega_j) = \frac{1}{(2 \pi)^D} \int\limits_\Upsilon d \mathbf{y^\prime} \int\limits_\Upsilon d \mathbf{y^{\prime\prime}} \ \text{exp}(- i \boldsymbol{\ell} \mathbf{y^\prime}) \Phi(\mathbf{y^\prime},\mathbf{y^{\prime\prime}}; \omega_j) \text{exp}(i \mathbf{k} \mathbf{y^{\prime\prime}}), \ D = 2, 3.
		\end{equation} 	
		
\end{enumerate}

To derive relations (\ref{Phi_int}), (\ref{f_int_new}), it is convenient to consider the above integral equations (\ref{F_F0_int})-(\ref{f_int}) in operator form. For example, the operator form of equation (\ref{F_F0_int}) can be written as
\begin{equation}
	\label{F_F0_operators}
	\hat{\text{G}}_0 (\hat{\mathscr{F}} - \hat{\mathscr{F}}_0) \hat{\text{G}} = \hat{\text{G}} - \hat{\text{G}}_0 ,
\end{equation} 	
\noindent
where each of operators is defined as 
$ (\hat{\mathscr{F}} - \hat{\mathscr{F}}_0) (\bullet) = 
\int\limits_\Upsilon
(\mathscr{F} - \mathscr{F}_0)(\mathbf{y^\prime},\mathbf{y^{\prime\prime}}; \omega_j) 
(\bullet)
d \mathbf{y^{\prime\prime}}
$,
$ \hat{\text{G}} = 
\int\limits_\Upsilon
G(\mathbf{y},\mathbf{x}; \omega_j) 
(\bullet)
d \mathbf{x}
$, 
$ \hat{\text{G}}_0 = 
\int\limits_\Upsilon
G_0(\mathbf{y} - \mathbf{x}; \omega_j) 
(\bullet)
d \mathbf{x}
$.
Then from (\ref{F_F0_operators}) directly follows the expression for the solution (\ref{F_F0_int}) in the form \cite{label30_1} $ (\hat{\mathscr{F}} - \hat{\mathscr{F}}_0) = \hat{\text{G}}^{-1}_0 - \hat{\text{G}}^{-1} $. 
By using operator expressions similar to (\ref{F_F0_operators}) for integrals (\ref{A_int}), (\ref{u_int}) and substituting the results in (\ref{f_int}), it is possible to obtain the following simplified expression for calculating the scattering amplitude:
\begin{equation}
\label{f_operators_solution}
\hat{f} = \hat{\text{U}}^\text{T}_\text{L} 
\left[ 
\hat{\text{G}}^{-1}_0 (\hat{\text{G}} - \hat{\text{G}}_0) \hat{\text{G}}^{-1}_0 
\right]
\hat{\text{U}}_\text{K},
\end{equation} 	
\noindent
where elements of matrices $\hat{f}$, $\hat{\text{U}}_\text{L}$, $\hat{\text{U}}_\text{K}$ are determined by the discretized values of functions $f(\mathbf{k}, \boldsymbol{\ell}; \omega_j)$, $\text{exp}(- i \boldsymbol{\ell} \mathbf{y})$, $\text{exp}(i \mathbf{k} \mathbf{y})$, respectively; the subscript "T" denotes transposition. 
Note, that the expression in square brackets in (\ref{f_operators_solution}) is the discretized analogue of the solution to the integral equation (\ref{Phi_int}). As a result, in numerical modeling, the scattering amplitude recalculation is reduced to one matrix equation (\ref{f_operators_solution}), which is implemented in the form of one line of program code.

The idea about the principal possibility of simplifying the recalculation procedure of $G_0(\mathbf{r}, \mathbf{x}; \omega_j)$, $G(\mathbf{r}, \mathbf{x}; \omega_j)$ into scattering amplitudes $f(\mathbf{k}, \boldsymbol{\ell}; \omega_j)$ was mentioned by O.D.~Rumyantseva during a discussion with author of this work, without specifying the details of a possible numerical implementation. 
The obtained relations (\ref{Phi_int}), (\ref{f_int_new}), (\ref{f_operators_solution}) and given below (\ref{Phi_double_angular_harm_new}), (\ref{f_double_angular_harm_new}) are not the only possible way for solving the considered problem (see, for example \cite{label30_2}).

Even greater efficiency in numerical modeling can be achieved when the contour $\Upsilon$ has a circle form. 
In this case, it is possible to avoid discretization of integral equations (\ref{Phi_int}), (\ref{f_int_new}) in the coordinate space and consider the problem in the space of angular harmonics \cite{label14, label23}. 
The numerical implementation of relations (\ref{Phi_int}), (\ref{f_int_new}) in the coordinate domain needs the adequate discretization of integrals, which requires additional control of accuracy for such discretization, involving, in general case, various interpolation algorithms that complicates the implementation of reconstruction procedure. This difficulty can be avoided by moving to the space of angular harmonics, which is Fourier conjugate to the space of angles that specify the positions of spatial points $\mathbf{x}, \mathbf{y}$ on the corresponding boundary $\Upsilon$.

Let the contour $\Upsilon$ has a circle shape with a radius $R_0$ with a center $\it O$; then in a polar coordinate system with the same center $\it O$ one can get: 
$\mathbf{x} = \left\{R_0, \phi_\mathbf{x}\right\}$, 
$\mathbf{y} = \left\{R_0, \phi_\mathbf{y}\right\}$, 
$\mathbf{y^{\prime}} = \left\{R_0, \phi_\mathbf{y^{\prime}}\right\}$, 
$\mathbf{y^{\prime\prime}} = \left\{R_0, \phi_\mathbf{y^{\prime\prime}}\right\}$. 
Below, the dependence on the parameter $R_0$ in the arguments of functions is omitted for brevity. 
For an arbitrary periodic function $g(\phi)$ with the period $2 \pi$, the transition to the angular spectrum $\tilde{g}(q)$ and the inverse transform are carried out by the formulas \cite{label14}:
\begin{equation}
\label{g_angular_harm}
\tilde{g}(q)=
\frac{1}{2 \pi} \int\limits_0^{2 \pi} g(\phi) \text{exp}(-i q \phi) d \phi, \
g(\phi) = \sum_{q = - \infty}^\infty \tilde{g}(q) \text{exp}(i q \phi), \
q \in \ \mathbb{Z},
\end{equation} 
\noindent
where $\mathbb{Z}$ is the set of integers; the sign "$\sim$" above the function means the Fourier transform with respect to the angular variable. 
For a function $g(\phi, \phi^\prime)$, depending on two angles $\phi$ and $\phi^\prime$, the double angular harmonic $\tilde{\tilde{g}}(q, q^\prime)$ is defined as
\begin{equation}
\label{g_double_angular_harm}
	\begin{split}
		& \tilde{\tilde{g}}(q, q^\prime) =
		\frac{1}{(2 \pi)^2} \int\limits_0^{2 \pi} \int\limits_0^{2 \pi}
		g(\phi, \phi^\prime) 
		\text{exp}(-i q \phi) \text{exp}(-i q^\prime \phi^\prime) 
		d \phi d \phi^\prime, \\
		& g(\phi, \phi^\prime) = 
		\sum_{q = - \infty}^\infty \sum_{q^\prime = - \infty}^\infty
		\text{exp}(i q \phi) \text{exp}(i q^\prime \phi^\prime).
	\end{split}
\end{equation} 
\noindent
Equation (\ref{Phi_int}) is rewritten in terms of angular harmonics as follows (similar to equation (31) from \cite{label14}):
\begin{equation}
\label{Phi_double_angular_harm_new}
\begin{split}
	& (2 \pi R_0)^2  
	\sum_{q^\prime_\mathbf{y} = - \infty}^\infty \ \sum_{q^{\prime\prime}_\mathbf{y} = - \infty}^\infty
	\tilde{\tilde{G}}_0(q_\mathbf{y}, q^\prime_\mathbf{y}; \omega_j)
	\tilde{\tilde{\Phi}}(-q^\prime_\mathbf{y}, - q^{\prime\prime}_\mathbf{y}; \omega_j)
	\tilde{\tilde{G}}_0(q^{\prime\prime}_\mathbf{y}, q_\mathbf{x}; \omega_j) = \\
	& = \tilde{\tilde{G}}(q_\mathbf{y}, q_\mathbf{x}; \omega_j) - 
	\tilde{\tilde{G}}_0(q_\mathbf{y}, q_\mathbf{x}; \omega_j); \ 
	q_\mathbf{x}, q_\mathbf{y} \in \mathbb{Z}; \ D = 2.
\end{split}
\end{equation} 
\noindent
It is convenient to solve the system (\ref{Phi_double_angular_harm_new}) with respect to $\tilde{\tilde{\Phi}}(-q^\prime_\mathbf{y}, - q^{\prime\prime}_\mathbf{y}; \omega_j)$. 
The double angular spectrum
$\tilde{\tilde{G}}_0(q_\mathbf{y}, q_\mathbf{x}; \omega_j)$ 
of classical Green's function \newline
$G_0(\mathbf{y}, \mathbf{x}; \omega_j) = - (i/4) \, H^{(1)}_0(k_{0j} |\mathbf{y} - \mathbf{x}|)$ 
of two-dimensional space, appearing in (\ref{g_double_angular_harm}), is calculated analytically and has no singularity, unlike 
$G_0(\mathbf{y}, \mathbf{x}; \omega_j)$ 
in coordinate space \cite{label14}: 
$\tilde{\tilde{G}}_0(q_\mathbf{y}, q_\mathbf{x}; \omega_j) = 
- (i/4) \, \delta_{q_\mathbf{x}, -q_\mathbf{y}} 
H^{(1)}_{q_\mathbf{y}}(k_{0j} R_0)
J_{q_\mathbf{y}}(k_{0j} R_0)
$, 
where where $\delta$ is the Konecker delta, and $H^{(1)}_{q_\mathbf{y}}$, $J_{q_\mathbf{y}}$ are the Hankel function of the first kind and the Bessel function both of $q_\mathbf{y}$-th order, respectively.  
It is also possible to get the analytical expression for the angular spectrum of plane wave field 
$u_0(\mathbf{y}, \mathbf{k}; \omega_j) = \text{exp}(i \mathbf{k} \mathbf{y})$ 
with a wave vector 
$\mathbf{k} = \left\{k_{0j}, \phi\right\}$ 
\cite{label14}: 
$ \tilde{u}_0(q_\mathbf{y}, \phi; \omega_j) = i^{q_\mathbf{y}} 
	J_{q_\mathbf{y}}(k_{0j} R_0) 
	\text{exp}(-i q_\mathbf{y} \phi)
	$. 
Finally, equation (\ref{f_int_new}), taking into account the representation 
$\mathbf{k} = \left\{k_{0j}, \phi\right\}$,
$\boldsymbol{\ell} = \left\{k_{0j}, \phi^\prime\right\}$, 
is rewritten in terms of angular variables as follows (similar to relation (40) from \cite{label14}):
\begin{equation}
\label{f_double_angular_harm_new}
	\begin{split}
		& f(\phi, \phi^\prime; \omega_j) = \\
		& = R_0^2 
		\sum_{q^\prime_\mathbf{y} = - \infty}^\infty \ \sum_{q^{\prime\prime}_\mathbf{y} = - \infty}^\infty
		\tilde{u}_0(q^\prime_\mathbf{y}, \phi^\prime + \pi; \omega_j)
		\tilde{\tilde{\Phi}}(-q^\prime_\mathbf{y}, - q^{\prime\prime}_\mathbf{y}; \omega_j)
		\tilde{u}_0(q^{\prime\prime}_\mathbf{y}, \phi; \omega_j).
	\end{split}
\end{equation} 
\noindent
The obtained relations (\ref{Phi_double_angular_harm_new}), (\ref{f_double_angular_harm_new}) solve the problem of recalculating the fields $G(\mathbf{r}, \mathbf{x}; \omega_j)$, measured at the boundary of tomography region, into scattering amplitudes \\ $f(\mathbf{k}, \boldsymbol{\ell}; \omega_j)$, which makes it possible to proceed to the description of main stages of the considered iterative algorithm \cite{label19}.

\renewcommand{\thesection}{\large 3}
\section{\large Main steps of acoustic scatterers reconstruction by using iterative algorithm of R.G.~Novikov}

The considered iterative algorithm \cite{label19} uses the well-known relationship between the scattering amplitude 
$f(\mathbf{k}, \boldsymbol{\ell}; \omega_j)$ 
and the inhomogeneity (scatterer) $v(\mathbf{r}, \omega_j)$:
\begin{equation}
\label{f_and_scatt_relation}
f(\mathbf{k}, \boldsymbol{\ell}; \omega_j) = 
\frac{1}{(2 \pi)^2} 
\int\limits_\mathfrak{R} \text{exp}(-i \boldsymbol{\ell} \mathbf{r})
v(\mathbf{r}, \omega_j)
u(\mathbf{r}, \mathbf{k}; \omega_j) 
d \mathbf{r},
\end{equation} 	
\noindent
where $u(\mathbf{r}, \mathbf{k}; \omega_j) $ is the field in the scattering domain 
$\mathfrak{R}$ 
(two-dimensional in the considered case) arising in response to the incident plane wave 
$u_0(\mathbf{r}, \mathbf{k}; \omega_j) = \text{exp}(i \mathbf{k} \mathbf{r})$. 
Relation (\ref{f_and_scatt_relation}) shows that the scattering amplitude 
$f(\mathbf{k}, \boldsymbol{\ell}; \omega_j)$ 
is determined in the form of a spatial Fourier transform of secondary sources 
$v(\mathbf{r}, \omega_j)
u(\mathbf{r}, \mathbf{k}; \omega_j)$ 
\cite{label2}, which arise in the scattering domain 
$\mathfrak{R}$. 
In the case when inhomogeneities are small (weak scatterers) and the multiple-scattering processes can be neglected, i.e. when the first Born approximation is valid, one can put 
$u(\mathbf{r}, \mathbf{k}; \omega_j) = 
u_0(\mathbf{r}, \mathbf{k}; \omega_j) = 
\text{exp}(i \mathbf{k} \mathbf{r})$ 
in (\ref{f_and_scatt_relation}), which leads to the so-called Born estimate of scatterer in the form:
\begin{equation}
\label{Born_approx}
\tilde{v}(\mathbf{k} - \boldsymbol{\ell}, \omega_j) =
f(\mathbf{k}, \boldsymbol{\ell}; \omega_j).
\end{equation}
\noindent
In this case, the scatterer 
$v(\mathbf{r}, \omega_j)$ 
is found from (\ref{Born_approx}) by the inverse Fourier transform:
\begin{equation}
\label{Fourier_trans_space_vectors}
v(\mathbf{r}, \omega_j) = 
\int\limits_{B_{2 k_{0j}}} 
\text{exp}(-i \boldsymbol{\xi} \mathbf{r})
\tilde{v}(\boldsymbol{\xi}, \omega_j)
d \boldsymbol{\xi}, \ \ 
\boldsymbol{\xi} = \mathbf{k} - \boldsymbol{\ell},
\end{equation}
\noindent
where the integration is carried out in the space of wave vectors 
$\boldsymbol{\xi}$ 
over the domain 
$B_{2 k_{0j}}$, 
which in the two-dimensional case is the circle with center at the origin of coordinates and with radius
$2 k_{0j}$. 
In the case when the inhomogeneity
$v(\mathbf{r}, \omega_j)$ 
significantly distorts the incident field and the first Born approximation is not valid, the scatterer estimate can be based on the following relation (see (2.10) in \cite{label19}):
\begin{equation}
\label{Iter_algo_main_ratio}
\begin{split}
\tilde{v}(\mathbf{k} - \boldsymbol{\ell}, \omega_j) = 
& f(\mathbf{k}, \boldsymbol{\ell}; \omega_j) - \\
& - \frac{1}{(2 \pi)^2} 
\int\limits_\mathfrak{R} 
\text{exp}(-i \boldsymbol{\ell} \mathbf{r})
v(\mathbf{r}, \omega_j)
\left[
u(\mathbf{r}, \mathbf{k}; \omega_j) -
\text{exp}(i \mathbf{k} \mathbf{r})
\right]
d \mathbf{r},
\end{split}
\end{equation}
\noindent
which, taking into account (\ref{f_and_scatt_relation}), is just the identity 
$\tilde{v}(\mathbf{k} - \boldsymbol{\ell}, \omega_j) = 
f(\mathbf{k}, \boldsymbol{\ell}; \omega_j) -
f(\mathbf{k}, \boldsymbol{\ell}; \omega_j) +
\tilde{v}(\mathbf{k} - \boldsymbol{\ell}, \omega_j) =
\tilde{v}(\mathbf{k} - \boldsymbol{\ell}, \omega_j)
$. Relation (\ref{Iter_algo_main_ratio}) allows one to propose a simple in terms of numerical implementation, but at the same time very effective in the sense of convergence \cite{label19}, the procedure for iterative estimation of scatterer
$
\tilde{v}(\mathbf{k} - \boldsymbol{\ell}, \omega_j)
$:  
\begin{equation}
\label{Iterative_algo}
\tilde{v}^{(n)}(\mathbf{k} - \boldsymbol{\ell}, \omega_j) = 
\tilde{v}^{(n-1)}(\mathbf{k} - \boldsymbol{\ell}, \omega_j) + 
f(\mathbf{k}, \boldsymbol{\ell}; \omega_j) - 
f^{(n-1)}(\mathbf{k}, \boldsymbol{\ell}; \omega_j),
\end{equation}
\noindent
where
$n = 1, 2, 3 ...$ 
is the number of iteration step, 
$f^{(n-1)}(\mathbf{k}, \boldsymbol{\ell}; \omega_j)$ 
is the scattering amplitude, which should be calculated at the $(n-1)$-th step by using the scatterer estimate
$
\tilde{v}^{(n-1)}(\mathbf{k} - \boldsymbol{\ell}, \omega_j)
$. The Born estimate (\ref{Born_approx}) can be chosen as an initial approximation
$
\tilde{v}^{(0)}(\mathbf{k} - \boldsymbol{\ell}, \omega_j) = 
f(\mathbf{k}, \boldsymbol{\ell}; \omega_j)
$, but this is not necessary; the choice of initial approximation can be based, for example, on a priori information about the functions being reconstructed.

Thus, the main steps of the considered iterative algorithm are as follows:

$\underline{\text{Step 1}}$. Recalculation of scattering data measured at the boundary of tomography area in the form of acoustic fields into scattering amplitudes
$
f(\mathbf{k}, \boldsymbol{\ell}; \omega_j)
$, based on (\ref{Phi_double_angular_harm_new}), (\ref{f_double_angular_harm_new}).

$\underline{\text{Step 2}}$. The choice of initial scatterer estimate 
$
\tilde{v}^{(0)}(\mathbf{k} - \boldsymbol{\ell}, \omega_j)
$, which can be the first Born approximation (\ref{Born_approx}).

$\underline{\text{Step 3}}$. The transition from the space of wave vectors to the coordinate space based on (\ref{Fourier_trans_space_vectors}), i.e. calculation 
$
v^{(n-1)}(\mathbf{r}, \omega_j)
$ 
from data
$
\tilde{v}^{(n-1)}(\mathbf{k} - \boldsymbol{\ell}, \omega_j)
$. At this step, it is required to make a Fourier transform of a function
$
\tilde{v}^{(n-1)}(\mathbf{k} - \boldsymbol{\ell}, \omega_j)
$
given on a non-uniform grid in the wavenumber space
$
\boldsymbol{\xi}
$. A similar problem has already been considered earlier in \cite{label31}.

At this step the implementation of (\ref{Fourier_trans_space_vectors}) can be made over a domain $B_{2 \tau k_{0j}}$ with a smaller radius $2 \tau k_{0j}$, where values of parameter $\tau \le 1$ can be changed during iterations to improve its convergence \cite{label19, label20}. This is equivalent to adaptive filtering of scatterer spatial spectrum aimed at a stepwise reconstruction of its high-frequency components during iterations.

$\underline{\text{Step 4}}$. Calculation of acoustic fields
$    
G^{(n-1)}(\mathbf{y}, \mathbf{x}; \omega_j)
$ at the boundary $\Upsilon$ in the presence of inhomogeneity
$
v^{(n-1)}(\mathbf{r}, \omega_j)
$
in the tomography area. This problem is solved on the basis of two Lippmann-Schwinger equations, one of which allows finding the fields
$
G^{(n-1)}(\mathbf{r}, \mathbf{x}; \omega_j)
$
at the internal points of scattering region
$
\mathbf{r} \in
\mathfrak{R}
$, and the second equation gives the fields
$
G^{(n-1)}(\mathbf{y}, \mathbf{x}; \omega_j)
$
at the boundary of investigated  region
$
\mathbf{x}, \mathbf{y} \in \Upsilon
$
by using
$
G^{(n-1)}(\mathbf{r}, \mathbf{x}; \omega_j)
$
\cite{label32}:
\begin{equation}
\label{Lippman_Schvinger_2equations}
\begin{split}
& G^{(n-1)}(\mathbf{r}, \mathbf{x}; \omega_j) = 
G_0(\mathbf{r}, \mathbf{x}; \omega_j) +
\int\limits_\mathfrak{R} 
G_0(\mathbf{r}, \mathbf{r^\prime}; \omega_j)
v^{(n-1)}(\mathbf{r^\prime}, \omega_j)
G^{(n-1)}(\mathbf{r^\prime}, \mathbf{x}; \omega_j)
d \mathbf{r^\prime}, \\
& G^{(n-1)}(\mathbf{y}, \mathbf{x}; \omega_j) = 
G_0(\mathbf{y}, \mathbf{x}; \omega_j) +
\int\limits_\mathfrak{R} 
G_0(\mathbf{y}, \mathbf{r^\prime}; \omega_j)
v^{(n-1)}(\mathbf{r^\prime}, \omega_j)
G^{(n-1)}(\mathbf{r^\prime}, \mathbf{x}; \omega_j)
d \mathbf{r^\prime}. 
\end{split}
\end{equation}
\noindent 
It should be noted that the solution to the direct problem of acoustic scattering, i.e. finding scattered fields for a known inhomogeneity, can be solved by any other methods \cite{label26, label27, label28}. The choice of method for solving the direct problem is not principal; it is only important that it allows calculating the fields in a reasonable time and takes into account the multiple-scattering processes with high accuracy.

$\underline{\text{Step 5}}$. Calculation of scattering amplitude
$
f^{(n-1)}(\mathbf{k}, \boldsymbol{\ell}; \omega_j)
$
from
$
G^{(n-1)}(\mathbf{y}, \mathbf{x}; \omega_j)
$,
$
G_0(\mathbf{y}, \mathbf{x}; \omega_j)
$
by using (\ref{Phi_double_angular_harm_new}), (\ref{f_double_angular_harm_new}).

$\underline{\text{Step 6}}$. Calculation of $ n $-th scatterer estimate
$
\tilde{v}^{(n)}(\mathbf{k} - \boldsymbol{\ell}, \omega_j)
$
from (\ref{Iterative_algo}), by using the functions
$
\tilde{v}^{(n-1)}(\mathbf{k} - \boldsymbol{\ell}, \omega_j)
$,
$
f^{(n-1)}(\mathbf{k}, \boldsymbol{\ell}; \omega_j)
$.

$\underline{\text{Step 7}}$. Repeating Steps 3-6 for the subsequent iterations $n = 2, 3, 4...$. The criterion for stopping the iterative process can be, for example, the achievement of a given value of discrepancy between the initial scattering amplitudes $f(\mathbf{k}, \boldsymbol{\ell}; \omega_j)$ and the estimate $f^{(n)}(\mathbf{k}, \boldsymbol{\ell}; \omega_j)$ obtained at the current iteration step.

A schematic representation of the described above iterative algorithm can be as follows:


\[
G(\mathbf{y}, \mathbf{x}; \omega_j) 
\xrightarrow[]{(\ref{Phi_double_angular_harm_new}), (\ref{f_double_angular_harm_new})}
f(\mathbf{k}, \boldsymbol{\ell}; \omega_j)
\xrightarrow[]{(\ref{Born_approx})}
\tilde{v}^{(0)}(\mathbf{k} - \boldsymbol{\ell}, \omega_j)
\xrightarrow[]{(\ref{Fourier_trans_space_vectors})} 
\]

\[
\xrightarrow[]{(\ref{Fourier_trans_space_vectors})}
v^{(n-1)}\tikzmark{b}(\mathbf{r}, \omega_j) 
\xrightarrow[]{(\ref{Lippman_Schvinger_2equations})}
G^{(n-1)}(\mathbf{y}, \mathbf{x}; \omega_j) 
\xrightarrow[]{(\ref{Phi_double_angular_harm_new}), (\ref{f_double_angular_harm_new})}
f^{(n-1)}(\mathbf{k}, \boldsymbol{\ell}; \omega_j)
\xrightarrow[]{(\ref{Iterative_algo})}
\tilde{v}^{(n)}\tikzmark{a}(\mathbf{k} - \boldsymbol{\ell}, \omega_j)
\tikz[overlay,remember picture]
{\path[draw,->,square arrow] (a.south) to node{{\footnotesize (18)}} (b.south) ;}
\]

\vspace{30pt}

It should be noted that there are many known iterative methods for solving inverse problems of acoustic scattering \cite{label2}. The considered algorithm \cite{label19} is distinguished, first of all, because for it there are a number of rigorously proven mathematical results. For example, mathematically rigorous convergence estimates were obtained even in the case of incomplete (fragmentary data). The Lipschitz stability is also proved, which is especially valuable for practical applications. At the same time, for other iterative algorithms known to the author of this work, the similar properties have been shown based on physical considerations, as well as on the results of individual physical or numerical experiments, but have not been proven rigorously. Another advantage of the considered iterative algorithm is the simplicity of its numerical implementation. Indeed, it is necessary just to make the Fourier transform (\ref{Fourier_trans_space_vectors}) and solve the direct problem, i.e. calculate fields that have passed through the known inhomogeneous medium. Finally, the considered algorithm seems to be promising for solving the so-called multichannel inverse problems, the acoustic analogue of which is the ocean nonadiabatic mode tomography, which takes into account the multichannel scattering of hydroacoustic modes \cite{label16}.

At the same time, mathematically rigorous results for the considered iterative algorithm \cite{label19} were obtained for inverse problems for the Schr{\"o}dinger equation. It is of interest to study numerically the possibilities of this approach in solving acoustic inverse problems for the Helmholtz equation.

\renewcommand{\thesection}{\large 4}
\section{\large Numerical modeling}

In the numerical simulation, a two-dimensional region $V_S$ of cylindrical shape with the radius $R_0$ surrounded by 60 receiving-emitting transducers was considered.
Parameters of models discussed below are mostly illustrative and are chosen to study capabilities of the considered iterative algorithm. The radius of region $V_S$ was assumed to be equal to $R_0 = 4 \lambda_{01}$, where $\lambda_{01}$ is the wavelength in the background environment corresponding to the lowest considered frequency $\omega_1$ and expressed in relative length sampling units (l.s.u.): $\lambda_{01}$~=~8~l.s.u.  
Acoustic fields, which are initial scattering data, were calculated by solving the Lippmann-Schwinger equations (\ref{Lippman_Schvinger_2equations}) in the presence of inhomogeneity $v(\mathbf{r}, \omega_j)$ in the tomography area. 
To characterize the accuracy of estimates $\hat{v}^{(n)}(\mathbf{r}, \omega_j)$ obtained at the $n$-th iteration step, the relative root-mean-square (rms) reconstruction errors (discrepancies for the solution) are calculated over the entire tomography region $V_S$: 
\[
\delta^{(n)}_v \equiv \sqrt{\int\limits_{V_S} \bigl| \hat{v}^{(n)}(\mathbf{r}, \omega_j) - v(\mathbf{r}, \omega_j) \bigr|^2 d \mathbf{r}} \bigg/ \sqrt{\int\limits_{V_S} \bigl| v(\mathbf{r}, \omega_j) \bigr|^2 d \mathbf{r}}.
\]

Reconstruction results obtained by using the above-mentioned functional-analytical algorithm \cite{label21, label21} and corresponding discrepancies will be denoted as $\hat{v}(\mathbf{r}, \omega_j)$ and $\delta_v$, i.e. without superscripts $ n $. 

To describe the strength of considered scatterers, i.e. to estimate how strongly they distort the incident acoustic field, the values of additional phase shifts are calculated
\[
\Delta \psi = k_{0j} \int\limits_{l_\mathfrak{R}} \frac{\Delta c(\mathbf{r}) \big/ c_0}{1 + \Delta c(\mathbf{r}) \big/ c_0} d l_\mathbf{r},
\]
\noindent
here $\Delta c(\mathbf{r}) \big/ c_0 \equiv \left\{ c(\mathbf{r}) - c_0 \right\} \big/ c_0$ is the relative sound speed contrast, $d l_\mathbf{r}$ is the length element of trajectory $l_\mathfrak{R}$ in the vicinity of point $\mathbf{r}$. The norm of scattering data is also calculated in the form
\[
\| f(\phi, \phi^\prime) \| \equiv \sqrt{\int\limits_0^{2 \pi} d \phi \int\limits_0^{2 \pi} d \phi^\prime \ \bigl| f(\phi, \phi^\prime) \bigr|^2 },
\]
\noindent
where $\phi$, $\phi^\prime$ are angular components of wave vectors $\mathbf{k}$, $\boldsymbol{\ell}$ (see (\ref{f_double_angular_harm_new}), dependence of scattering amplitude on $\omega_j$ here and below is omitted). The norm of scattering data, together with the additional phase shift, characterizes the scatterer strength. This norm was used in \cite{label10} to determine the sufficient condition of convergence for the functional-analytical algorithm as $\| f(\phi, \phi^\prime) \| < 1 \big/ (3 \pi)$. However, later it was shown (see \cite{label23} and references therein) that the functional-analytical algorithm remains stable even when the norm of data $\| f(\phi, \phi^\prime) \|$ is an order of magnitude or more higher than the threshold value $1 / (3 \pi)$. 

\begin{figure}[t!]
	\centerline{\epsfig{file=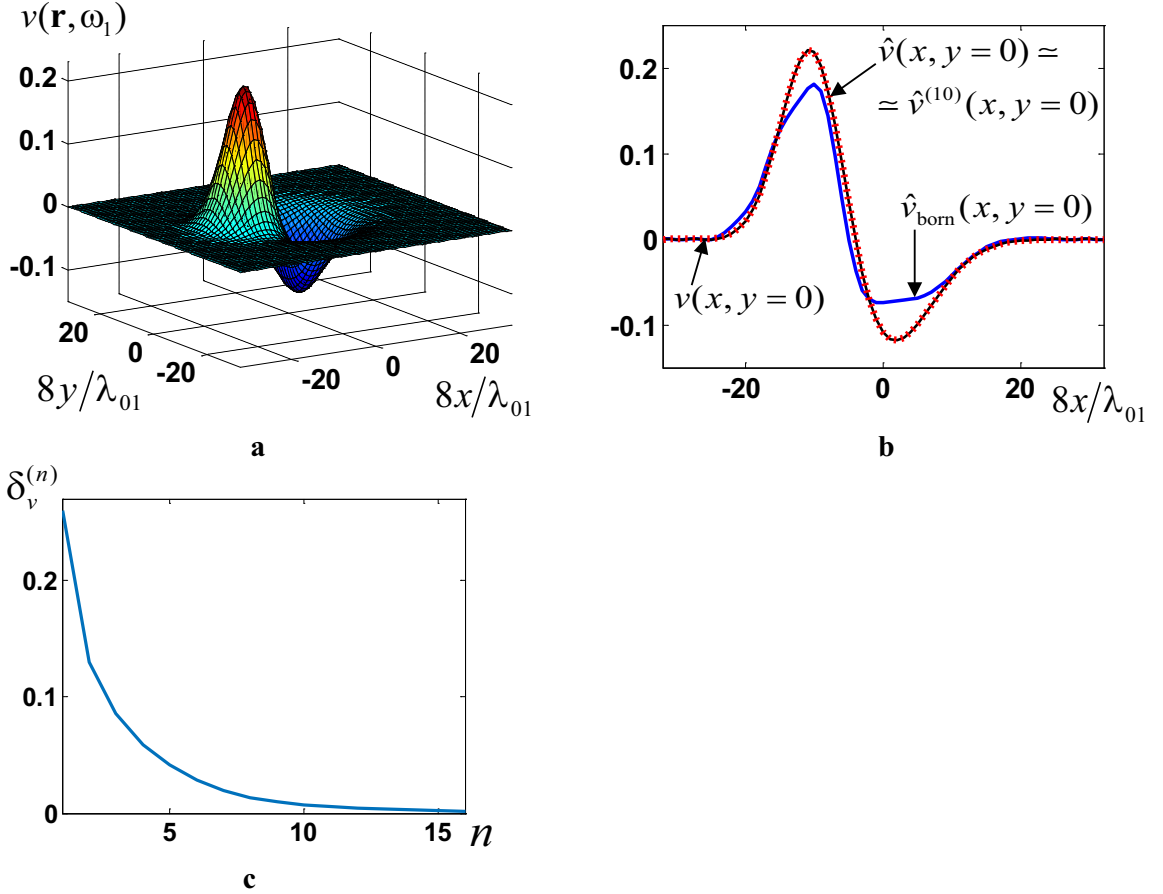, width = 15.5 cm}}
	\caption{ General view of true scatterer $v$ (a), for which the relative contrast of sound speed $\Delta c(\mathbf{r}) \big/ c_0$ ranges from -0.084 to 0.25, maximum additional phase shift is $\Delta \psi \approx 0.46 \pi$, norm of scattering data is $\| f(\phi, \phi^\prime) \| \approx 11 / (3 \pi)$, dimensionless coefficient is $\text{A}_0 = 0.43$; 
		\protect\\
		- central cross sections of true scatterer $ v $ (b, thin solid line), reconstruction results obtained by the iterative and functional-analytical algorithms, which are visually identical $ \hat{v} \simeq \hat{v}^{(10)} $ (b, dotted line), and the Born estimate $ \hat{v}_{\text{born}}$  (b, thick solid line);
		\protect\\
		- dependence of discrepancy for the solution $ \delta_v^{(n)} $ on the iteration number $ n $ (c).
	}
	\label{figure_2}
\end{figure}

\begin{figure}[t!]
	\centerline{\epsfig{file=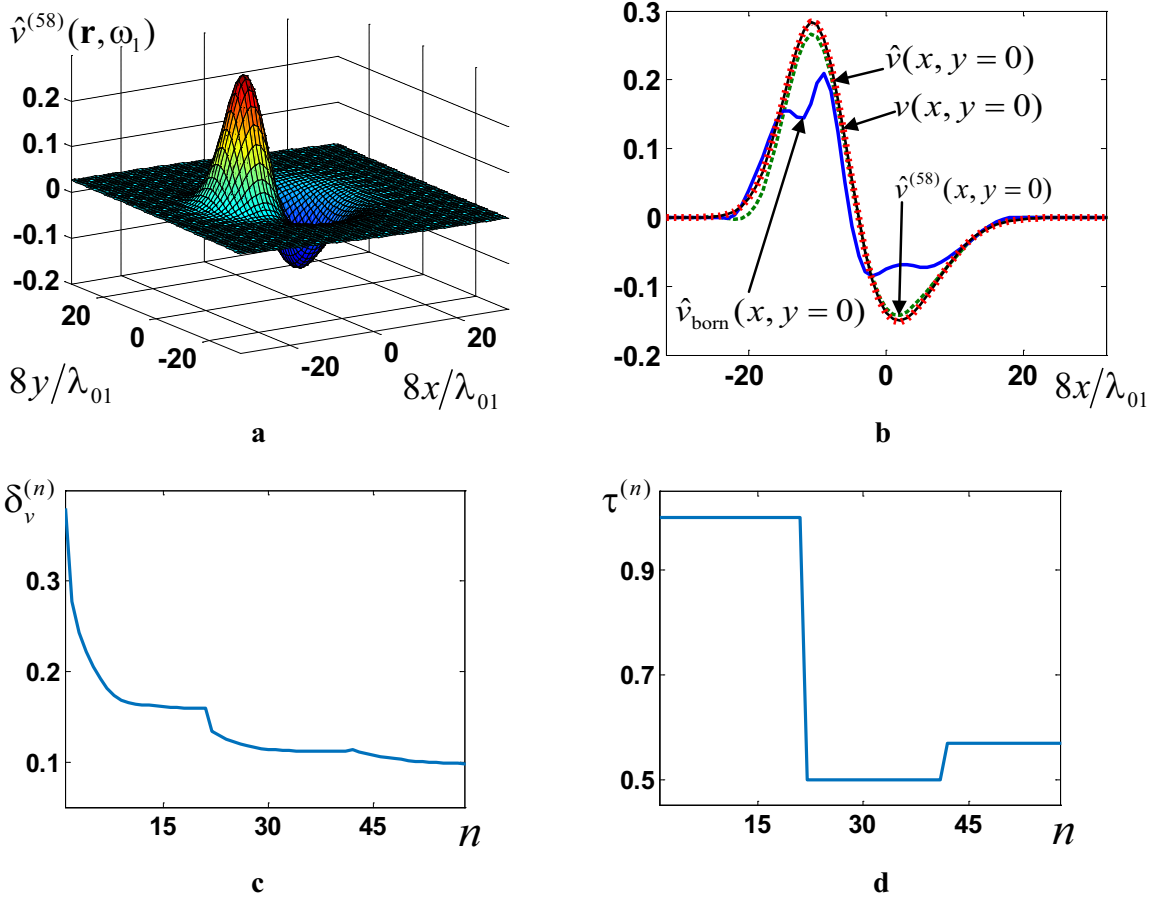, width = 15.5 cm}}
	\caption{ General view of the scatterer $ \hat{v}^{(58)}$ reconstructed at 58-th iteration (а), for which the relative contrast of sound speed $\Delta c(\mathbf{r}) \big/ c_0$ ranges from -0.1 to 0.36, maximum additional phase shift is $\Delta \psi \approx 0.6 \pi$, norm of scattering data is $\| f(\phi, \phi^\prime) \| \approx 13.4 / (3 \pi)$, dimensionless coefficient is $\text{A}_0 = 0.55$; 
		\protect\\
		- central cross sections of true scatterer $ v $ (b, thin solid line), reconstruction results obtained by the functional-analytical algorithm $ \hat{v}$ (b, dotted line), by the iteration method $ \hat{v}^{(58)}$ (b, dash line) and  the Born estimate  $ \hat{v}_{\text{born}}$  (b, thick solid line);
		\protect\\
		- dependence of discrepancy $ \delta_v^{(n)} $ (c) and the parameter of filtration $ \tau^{(n)} $ (d) on the iteration number $ n $.}
	\label{figure_3}
\end{figure}

The scatterers $v(\mathbf{r}, \omega_j)$ shown in figures~\ref{figure_2},~\ref{figure_3} consist of two Gaussian-shape inhomogeneities with different amplitude values, different sizes and located at different distances from the center of tomographic region:
\begin{equation}
\label{scatt_rec_1}
v(\mathbf{r}, \omega_j) = \text{A}_0 \ k^2_{0j}
\left[ 
\text{exp}(-|\mathbf{r} - \mathbf{r}^\prime| / \sigma^{\prime \, 2}) -   
0.5 \ \text{exp}(-|\mathbf{r} - \mathbf{r}^{\prime\prime}| / \sigma^{\prime\prime \, 2})
\right],
\end{equation}
\noindent
here $\mathbf{r}^\prime = \left\{ -10 / 8 \lambda_{01}, 0 \right\}$, $\mathbf{r}^{\prime\prime} = \left\{ 0, 0 \right\}$, $\sigma^\prime = 2/3 \sigma^{\prime\prime} = 0.8 \lambda_{01}$. To change the strength of scatterer defined in (\ref{scatt_rec_1}), various values of dimensionless coefficient $\text{A}_0$ are used below, while other parameters in (\ref{scatt_rec_1}) did not change. For the considered class of scatterers, the maximum additional phase shift $\Delta \psi$ occurs when a wave propagates along the $OX$ axis through the central cross section of reconstructed inhomogeneities; this phase shift will be calculated further.

Figure~\ref{figure_2} shows results of scatterer reconstruction with coefficient 
$
\text{A}_0 = 0.43
$. 
In this case, the additional phase shift is
$
\Delta \psi \approx 0.46 \pi
$
in the cross section ($y = 0$) with a positive velocity contrast
$
\Delta c(\mathbf{r}) \big/ c_0 > 0
$
and
$
\Delta \psi \approx -0.26 \pi
$
in the cross section with a negative velocity contrast
$
\Delta c(\mathbf{r}) \big/ c_0 < 0
$
(figure~\ref{figure_2}).  
Thus, the considered scatterer is strong enough and cannot be reconstructed with acceptable accuracy within the Born approximation. 
Figure~\ref{figure_2}b shows the cross section of the function
$
\hat{v}^{(10)}(\mathbf{r}, \omega_1)
$, reconstructed by the discussed iterative algorithm in $n$~=~10 steps at a single frequency $\omega_j = \omega_1$, the cross section of estimate
$
\hat{v}(\mathbf{r}, \omega_1)
$, obtained by the functional-analytical algorithm, as well as the cross section of true distribution
$
v(\mathbf{r}, \omega_1)
$. 
As it can be seen in figure~\ref{figure_2}b, the location of inhomogeneities, their shape and amplitudes are reconstructed with high accuracy in the absence of noises in the initial data. 
The values of discrepancies
$
\delta_v \simeq \delta_v^{(10)} \approx 0.008
$
coincide for the functional-analytical approach and the iterative algorithm with an accuracy of hundredths of a percent. 
In this case, the norm of scattering data is
$
\| f(\phi, \phi^\prime) \| \approx 11 / (3 \pi)
$, which significantly exceeds the above-mentioned threshold value 
$1 / (3 \pi)$. 
It should be noted that an attempt to reconstruct the same scatterer in the first Born approximation gives unsatisfactory results, and the residual value becomes equal to
$
\delta_v^{(0)} \approx 0.26
$. 
As an illustration of iterations convergence, figure~\ref{figure_1}c shows the dependence of discrepancy $\delta_v^{(n)}$ on the iteration number $ n $, which shows that in the absence of noises in the initial data, only 3-5 iteration steps are required to decrease the residual value $\delta_v^{(n)}$ less than 0.05, when results of reconstruction and the true scatterer are practically indistinguishable.

\begin{figure}[ht!]
	\centerline{\epsfig{file=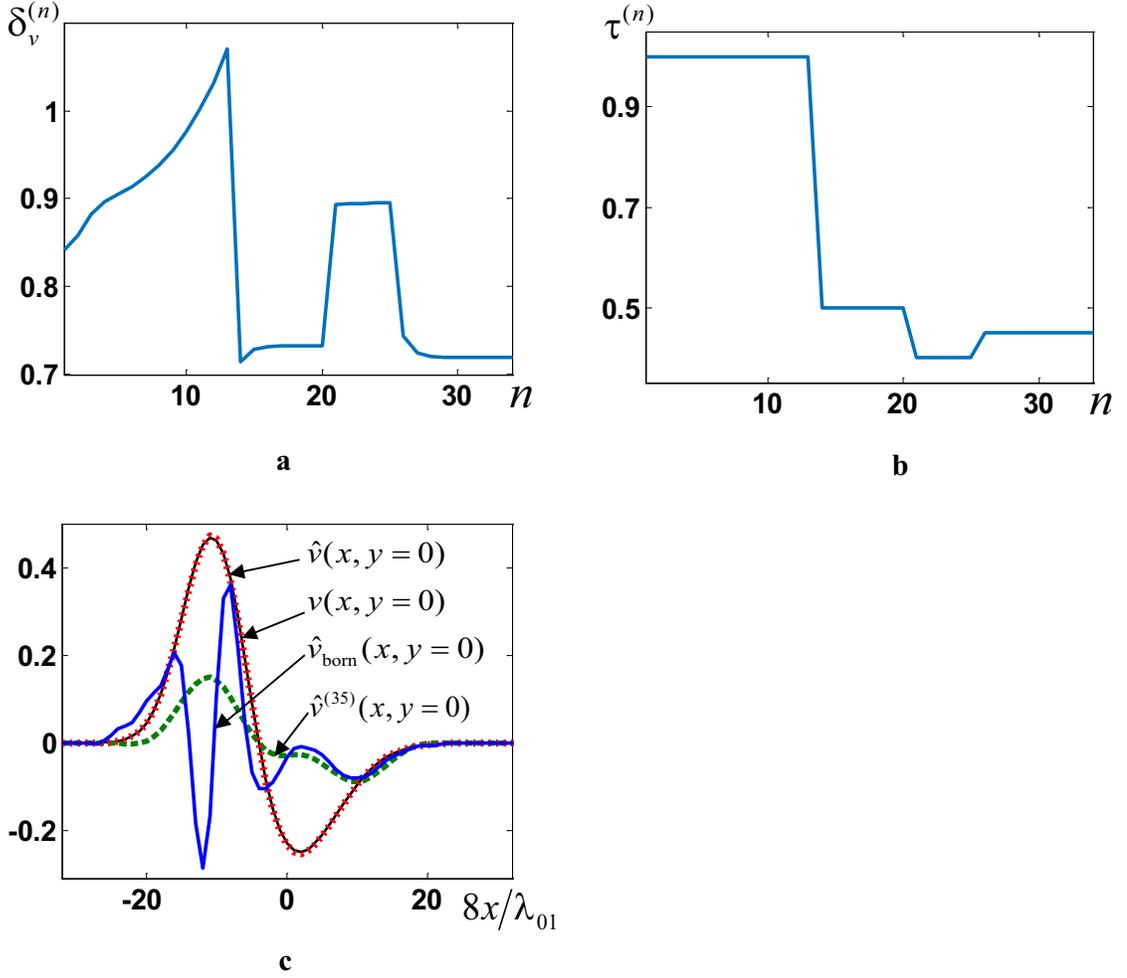, width = 15.5 cm}}
	\caption{ Dependence of discrepancy $ \delta_v^{(n)} $ (a) and parameter of filtration $ \tau^{(n)} $ (b) on the iteration number $ n $, when reconstructing a strong scatterer, for which values $\Delta c(\mathbf{r}) \big/ c_0$ ranges from - 0.16 to 1.08, maximum additional phase shift is $\Delta \psi \approx 1.1 \pi$, norm of scattering data is $\| f(\phi, \phi^\prime) \| \approx 19.3 / (3 \pi)$, dimensionless coefficient is $\text{A}_0 = 0.91$. After $ n $ = 35 iterations discrepancy is $ \delta_v^{(35)} \approx 0.72$, that is smaller than discrepancy of the Born estimate $ \delta_v^{(0)} \approx 0.84$, but significantly more than the functional-analytical algorithm result $ \delta_v \approx 0.012$; this shows the limitations of the iterative algorithm for the reconstruction of strong scatterers. The central cross sections of true scatterer $ v $ and its estimates are shown in (c): thin solid line shows $ v $, dotted line – $ \hat{v}$, dash line – $ \hat{v}^{(35)} $, thick solid line – $ \hat{v}_{\text{born}}$.  
	}
	\label{figure_4}
\end{figure}

With a further increase of scatterer strength, the results of reconstruction by the iterative algorithm deteriorate; however, the functional-analytical algorithm still yields a reconstruction with high accuracy. For the coefficient 
$
\text{A}_0 = 0.55
$ 
in (\ref{scatt_rec_1})
the additional phase shift is 
$
\Delta \psi \approx 0.6 \pi
$
in the cross section with a positive velocity contrast
$
\Delta c(\mathbf{r}) \big/ c_0 > 0
$
and
$
\Delta \psi \approx -0.33 \pi
$,
if
$
\Delta c(\mathbf{r}) \big/ c_0 < 0
$,
the norm of scattering data is
$
\| f(\phi, \phi^\prime) \| \approx 13.4 / (3 \pi)
$.
In this case the iterations begin to converge to a solution with a residual
$
\delta_v^{(20)} \approx 0.16
$
(see figure~\ref{figure_3}), which significantly exceeds the values obtained earlier for the case
$
\text{A}_0 = 0.43
$
(see figure~\ref{figure_2}). 
For the value
$
\text{A}_0 = 0.91
$
(additional phase shift is
$
\Delta \psi \approx 1.1 \pi
$,
if
$
\Delta c(\mathbf{r}) \big/ c_0 > 0
$,
and
$
\Delta \psi \approx -0.54 \pi
$,
if
$
\Delta c(\mathbf{r}) \big/ c_0 < 0
$,
the norm of scattering data is
$
\| f(\phi, \phi^\prime) \| \approx 19.3 / (3 \pi)
$
), the iterations even diverge (figure~\ref{figure_4}). 
One of the reason may be the accumulation of numerical errors, the magnitude of which becomes comparable to the contribution from each subsequent iterative addition. Such errors can arise, first of all, when solving a direct problem that requires a very small sampling step to take into account the multiple scattering processes on a small-scale but high-contrast elements of scatterers 
$
v^{(n-1)}(\mathbf{r}, \omega_j)
$,
which are reconstructed at intermediate iteration steps.
It should be noted that both in the case
$
\text{A}_0 = 0.55
$
and in the case
$
\text{A}_0 = 0.91
$,
the functional-analytical algorithm gives almost ideal results of reconstruction with discrepancies
$
\delta_v \approx 0.008
$
and
$
\delta_v \approx 0.012
$,
respectively; the values of these residuals characterizes the accuracy of solving the direct problem.
The possibilities of iterative algorithms in the reconstruction of middle strength scatterers and the difficulties encountered in the iterative reconstruction of strong scatterers were also discussed in \cite{label2, label32}.

As it was mentioned above in Section 3, convergence of iterations can be improved by introducing the additional filtering of spatial spectrum of scatterers \linebreak
$
\tilde{v}^{(n-1)}(\boldsymbol{\xi}, \omega_j)
$. 
Filtration of this kind is equivalent to considering in (\ref{Fourier_trans_space_vectors}) instead of the region $B_{2 k_{0j}}$ a circle with a smaller radius $2 \tau k_{0j}$, where $\tau \le 1$ \cite{label19}. 
Figure~\ref{figure_3} shows the reconstruction results obtained for the scatterer with coefficient
$
\text{A}_0 = 0.55
$,
when values of the parameter
$\tau^{(n)} \le 1$
were changed during iterations. 
The criterion for changing the value $\tau^{(n)}$ can be either the divergence of iterations, accompanied by an increase of the discrepancy $\delta^{(n)}_f$ for scattering data:
\[
\delta^{(n)}_f \equiv 
\sqrt{
	\int\limits_0^{2 \pi} d \phi \int\limits_0^{2 \pi} d \phi^\prime \ 
	\bigl| \hat{f}^{(n)}(\phi, \phi^\prime) - f(\phi, \phi^\prime) \bigr|^2 } 
\bigg/ 
\sqrt{
	\int\limits_0^{2 \pi} d \phi \int\limits_0^{2 \pi} d \phi^\prime \ 
	\bigl| f(\phi, \phi^\prime) \bigr|^2 } ,
\]
\noindent
here $\hat{f}^{(n)}(\phi, \phi^\prime)$ is the estimation of scattering amplitude at $ n $-th iteration step, or the convergence of iterations, when the value $\delta^{(n)}_f$ does not change significantly for several iterative steps. 
Another way of possible improvements of iteration convergence is the similar filtration of scattering amplitude inside a circle with radius $2 \tau k_{0j}$. This approach, in some sense, is equivalent to the stepwise inclusion method \cite{label32}. 
As it can be seen in figure~\ref{figure_3}c, when
$
\text{A}_0 = 0.55
$
the considered approach makes it possible to achieve the convergence of iterations with discrepancy
$
\delta^{(58)}_v \approx 0.098
$
that is smaller than the initial value 
$
\delta^{(0)}_v \approx 0.38
$, obtained in the Born approximation,
by more than 3.5 times . 
The final estimate
$
\hat{v}^{(58)}(\mathbf{r}, \omega_1)
$
is comparable in terms of reconstruction quality with the results of the functional-analytical method (see figure~\ref{figure_3}c).
In the case
$
\text{A}_0 = 0.55
$
the iterations converge (figure~\ref{figure_4}), but the final discrepancy turns out to be significantly larger than the reconstruction error obtained by the functional-analytical algorithm. 
The presented examples of numerical modeling correspond to the previously obtained theoretical estimates of accuracy of the iterative \cite{label19} and functional-analytical \cite{label21} approaches when solving the inverse problem for the Schr{\"o}dinger equation. In accordance with these estimates, the functional-analytical algorithm should give lower values of discrepancy for the solution in comparison with the iterative algorithm with increasing energy, which in acoustic case, in some sense, corresponds to an increase in the squared wavenumber. 
Thus, the obtained results of numerical simulation confirm for acoustic applications the validity of known theoretical estimates, thereby demonstrating the advantages of a rigorous functional-analytical algorithm in comparison with the iterative approach for recovering strong scatterers.

To analyze the resolution of the iterative algorithm, components with small sizes in the coordinate domain were added to the original scatterer (\ref{scatt_rec_1}):
\begin{equation}
\label{scatt_rec_2}
\begin{split}
v(\mathbf{r}, \omega_j) = & \text{A}_0 \ k^2_{0j}
[ 
\text{exp}(-|\mathbf{r} - \mathbf{r}^\prime| / \sigma^{\prime \, 2}) -   
0.5 \ \text{exp}(-|\mathbf{r} - \mathbf{r}^{\prime\prime}| / \sigma^{\prime\prime \, 2}) - \\
& - 0.5 \ \text{exp}(-|\mathbf{r} - \mathbf{r}^{\prime\prime\prime}| / \sigma^{\prime\prime\prime \, 2}) +
0.5 \ \text{exp}(-|\mathbf{r} - \mathbf{r}^{\prime\prime\prime\prime}| / \sigma^{\prime\prime\prime\prime \, 2})
],
\end{split}
\end{equation}
\noindent
where $\mathbf{r}^{\prime\prime\prime} = \left\{ -11 / 8 \lambda_{01}, 0 \right\}$, $\mathbf{r}^{\prime\prime\prime\prime} = \left\{ 1 / 4 \lambda_{01}, 0 \right\}$, $\sigma^{\prime\prime\prime} = \sigma^{\prime\prime\prime\prime} = 0.4 \lambda_{01}$, $\text{A}_0 = 1.1$. In this case, the scatterer contains components with a characteristic spatial scales about a quarter of a wavelength (figure~\ref{figure_5}), which is close to the limiting resolution of wave monochromatic methods for solving inverse scattering problems \cite{label23}. 
The considered scatterer (\ref{scatt_rec_2}) is quite strong: the additional phase shift is
$
\Delta \psi \approx 0.46 \pi
$
on a segment with a positive sound speed contrast
$
\Delta c(\mathbf{r}) \big/ c_0 > 0
$,
and
$
\Delta \psi \approx -0.19 \pi
$,
if
$
\Delta c(\mathbf{r}) \big/ c_0 < 0
$,
the norm of scattering data is
$
\| f(\phi, \phi^\prime) \| \approx 12.8 / (3 \pi)
$.
Despite this, it is possible to obtain acceptable reconstruction results with
$
\delta^{(17)}_v \approx 0.047
$
in $ n $~=~17 iteration steps without using the considered filtering in the space of wave vectors (see figure~\ref{figure_5}а). Thus, the considered iterative and functional-analytical algorithms have comparable resolution, while the functional-analytical approach gives a solution with a smaller discrepancy
$
\delta_v \approx 0.018
$.

Simultaneously with the resolution, the noise immunity of the discussed iterative procedure was analyzed. To make such investigation, a normally distributed random noise interference
$
n(\mathbf{y}, \mathbf{x}; \omega_j)
$, 
uncorrelated for the different directions of emitting and receiving transducers, with zero mean and with the rms amplitude deviation 
$
\sigma_{\text{ns}}(\omega_j) = 0.15 \bar{G}_{\text{sc}}(\omega_j)
$,
was introduced into the scattered fields
$
G_{\text{sc}}(\mathbf{y}, \mathbf{x}; \omega_j) \equiv G(\mathbf{y}, \mathbf{x}; \omega_j) - G_0(\mathbf{y}, \mathbf{x}; \omega_j)
$
separately for the real and imaginary parts and independently at different frequencies $\omega_j$.
The rms value $\bar{G}_{\text{sc}}$ of scattered data is defined as 
$
\bar{G}_{\text{sc}} \equiv 
\sqrt{
	\int\limits_{\Upsilon} d \mathbf{x} \int\limits_{\Upsilon} d \mathbf{y} \ 
	\bigl| G_{\text{sc}}(\mathbf{y}, \mathbf{x}; \omega_j) \bigr|^2 } 
\bigg/ 
\sqrt{
	\int\limits_{\Upsilon} d \mathbf{x} \int\limits_{\Upsilon} d \mathbf{y} } 
$,
and the input noise-to-signal amplitude ratio 
$
N / S \equiv 
\sqrt{
	\int\limits_{\Upsilon} d \mathbf{x} \int\limits_{\Upsilon} d \mathbf{y} \ 
	\bigl| n(\mathbf{y}, \mathbf{x}; \omega_j) \bigr|^2 } 
\bigg/ 
\sqrt{
	\int\limits_{\Upsilon} d \mathbf{x} \int\limits_{\Upsilon} d \mathbf{y} 
	\bigl| G_{\text{sc}}(\mathbf{y}, \mathbf{x}; \omega_j) \bigr|^2} 
$,
is
$
N / S \approx 0.21
$,
which even exceeds the level of interference expected, for example, in medical tomography experiments \cite{label33}. 
Nevertheless, the considered scatterer is reconstructed with acceptable accuracy using scattering data at one frequency $\omega_1$, corresponding to the wavelength $\lambda_{01}$ (see figure~\ref{figure_5}). 
The use of a multifrequency (impulse) sounding regime improves the results of reconstruction from noisy data \cite{label24}. In the considered case of scatterers, which are described by a perturbation of sound speed only, the multifrequency estimate at each iteration step is obtained by the simple averaging of functions
$
\hat{v}^{(n)}(\mathbf{r}, \omega_j) / \omega^2_j
$
over the frequencies $\omega_j$. As it follows from (\ref{Scatt_func}) these functions $\hat{v}^{(n)}(\mathbf{r}, \omega_j) / \omega^2_j$ do not depend on frequency. 
Due to the fact that at each frequency $ \omega_j$ the scattered fields $G_{\text{sc}}(\mathbf{y}, \mathbf{x}; \omega_j)$ contain independent realizations of noises $n(\mathbf{y}, \mathbf{x}; \omega_j)$, the simple averaging of functions $\hat{v}^{(n)}(\mathbf{r}, \omega_j) / \omega^2_j$ makes it possible to obtain the desired improvement in the noise-to-signal amplitude ratio $N/S$. 
For example, if $N/S \approx 0.21$, then the use of scattering data at 40 frequencies (in a real experiment, the amount of data can be even grater \cite{label33}) allows one to obtain reconstruction results that are indistinguishable from the results presented in figure~\ref{figure_5}b. 

\renewcommand{\thesection}{\large 5}
\section{\large Conclusions}

\begin{figure}[t!]
	\centerline{\epsfig{file=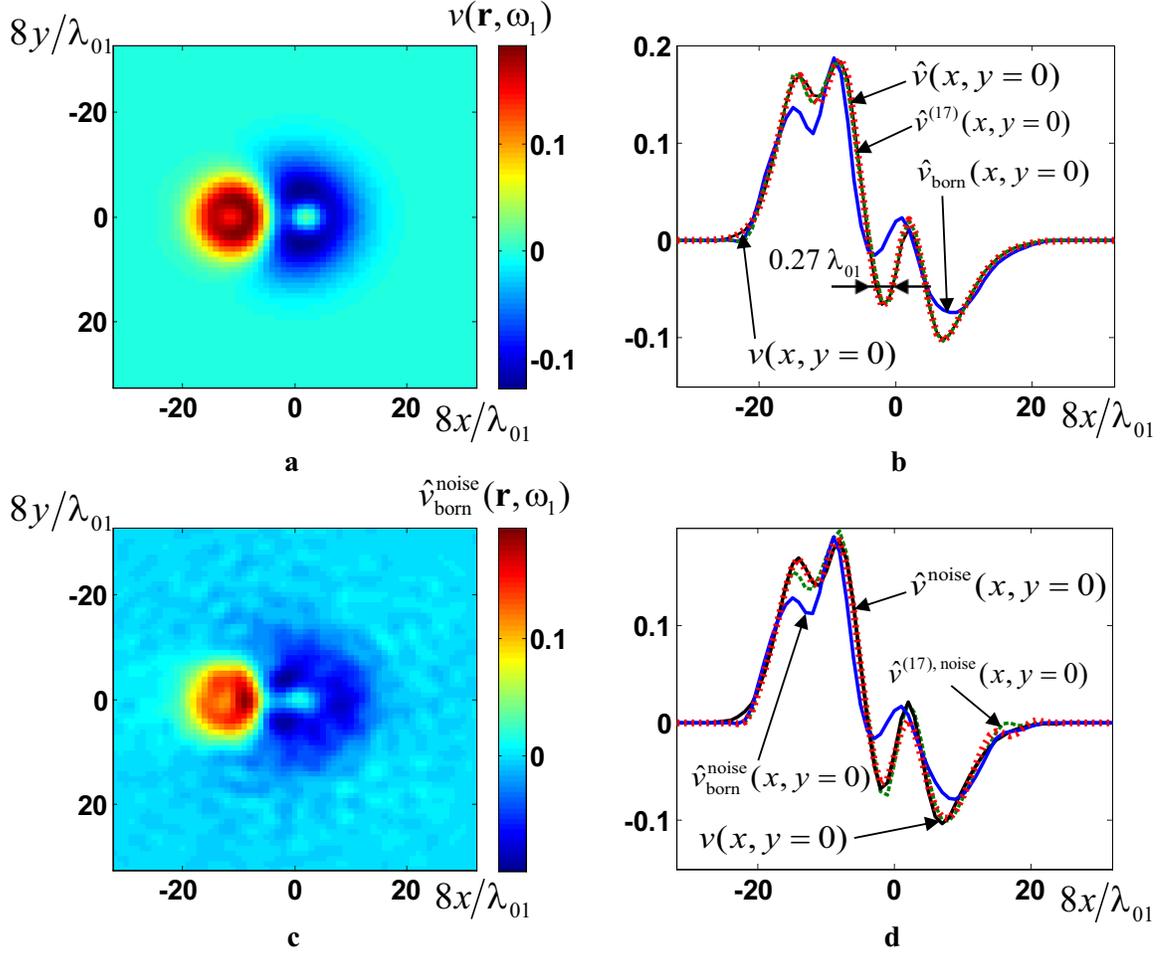, width = 15.5 cm}}
	\caption{ General view of middle strength scatterer with small size elements (а), for which the relative contrast of sound speed $\Delta c(\mathbf{r}) \big/ c_0$ ranges from - 0.07 to 0.19, maximum additional phase shift is $\Delta \psi \approx 0.46 \pi$, norm of scattering data is $\| f(\phi, \phi^\prime) \| \approx 12.8 / (3 \pi)$, dimensionless coefficient is $\text{A}_0 = 1.1$; 
		\protect\\
		- central cross sections of true scatterer $ v $ and its estimates obtained by using scattering data without noise are shown in (b): thin solid line shows $ v $, dotted line – $ \hat{v}$ (the discrepancy is $ \delta_v \approx 0.018$), dash line – $ \hat{v}^{(17)} $ ($ \delta_v^{(17)} \approx 0.047$), thick solid line – $ \hat{v}_{\text{born}}$ ($ \delta_v^{(0)} \approx 0.31$);
		\protect\\
		- general view of the Born estimate $ \hat{v}^\text{noise}_{\text{born}}(\mathbf{r}, \omega_1)$ (c) obtained by using noisy data at one frequency with noise rms amplitude deviation $\sigma_{\text{ns}}(\omega_j) = 0.15 \bar{G}_{\text{sc}}(\omega_j)$;
		\protect\\
		- central cross sections of true scatterer $ v $ and its estimates obtained by using noisy data are shown in (d): thin solid line shows $ v $, dotted line – $ \hat{v}^\text{noise}$ (the discrepancy is $ \delta_v \approx 0.018$), dash line – $ \hat{v}^{(17), \, \text{noise}} $ ($ \delta_v^{(17)} \approx 0.047$), thick solid line – $ \hat{v}^\text{noise}_{\text{born}}$ ($ \delta_v^{(0)} \approx 0.31$).
	}
	\label{figure_5}
\end{figure}

Results of this work show that the considered iterative algorithm \cite{label19} makes it possible to reconstruct acoustic scatterers of middle strength with accuracy, resolution and noise immunity comparable to reconstruction results of the functional-analytical approach \cite{label21, label22} (see figures~\ref{figure_2},~\ref{figure_5}). 
When recovering strong scatterers, a divergence of iterations is observed. 
To improve convergence of the considered iterative scheme, filtering of intermediate reconstruction results in the space of wave vectors was applied (figure~\ref{figure_3}). 
It should be noted that the scatterers, for which the divergence of iterations was observed (the maximum additional phase shift is
$
\Delta \psi \approx 1.1 \pi
$,
the norm of scattering data is
$
\| f(\phi, \phi^\prime) \| \approx 19.3 / (3 \pi)
$),
could be reconstructed by using the functional-analytical algorithm with high accuracy (figure~\ref{figure_4}) that corresponds to the known theoretical estimates of these algorithms convergence \cite{label19, label21}. 
Indeed, the main advantage of the functional-analytical algorithm is its mathematical rigor, which makes it possible to take into account the processes of multiple scattering when solving the inverse problem, thereby providing a non-iterative reconstruction of scatterers of different strengths with high accuracy. 
As a result, for the case of two-dimensional scalar Helmholtz equation, the efficiency of the iterative algorithm in reconstruction of middle strength scatterers and the advantages of the functional-analytical approach in reconstruction of strong scatterers with parameters close to acoustic tomography problems were demonstrated. 
A comparative numerical study of iterative and functional-analytical reconstruction algorithms is carried out in this work for the first time. 
The obtained results indicate perspectives of using the considered iterative algorithm for developing new acoustic tomography schemes that are flexible enough for use in various applications. The main advantage of the iterative approach in comparison with the functional-analytical algorithm is the ability to make reconstruction from incomplete scattering data \cite{label19}. This is especially important in three-dimensional inverse problems, an example of which is the nonadiabatic mode tomography of ocean \cite{label16}. 
The analysis of applicability of the iterative algorithm in such problems refers to prospects for further research.

\section{\large Acknowledgement}
The reported study was funded by RFBR and CNRS, project number 20-51-15004.



\begin{thebibliography}{}
\thispagestyle{myheadings}
\vspace*{-10mm}
\footnotesize

\bibitem{label1} {\it Physical principles of medical ultrasonics. Second edition.} Eds. C.R.~Hill , J.C.~Bamber , G.R.~ter Haar, London: John Wiley \&  Sons, Ltd, 2004.

\bibitem{label2} {V.A. Burov, O.D. Rumyantseva}, {\it Inverse Wave Problems of Acoustic Tomography. Part II: Inverse Problems of Acoustic Scattering} (in Russian), Lenand/URSS, Moscow, 2020.

\bibitem{label3} {W.~Munk, P.~Worcester, C.~Wunsch}, {\it Ocean Acoustic Tomography}, New York: Cambridge University Press, 1995.

\bibitem{label4} {V.V. Goncharov, V.Yu. Zaytsev, V.M.~Kurtepov, A.G.~Nechaev , A.I.~Khilko}, {\it Ocean Acoustic Tomography} (in Russian), IPF RAN, Nizhni Novgorod, 1997.

\bibitem{label5} {T.B.~Yanovskaya}, {\it Surface and Wave Tomography in Seismic Researches} (in Russian), Nauka, St. Petersburg, 2015.

\bibitem{label6} {A.L.~Sobisevich, D.A.~Presnov, A.S.~Shurup}, {\it Fundamentals of Passive Seismohydroacoustic Methods for Arctic Shelf Investigation}, Acoust. Phys. 67 (2021), no. 1, 62-82.

\bibitem{label7} {L.D.~Faddeev}, {\it Inverse problem of quantum scattering theory II}, J. of Soviet Math., 5 (1976), 334–396.

\bibitem{label8} {R.G.~Novikov}, {\it Construction of two-dimensional Schr{\"o}dinger operator with given scattering amplitude at fixed energy}, Theoret. and Math. Phys., 66 (1986), no. 2, 154–158.

\bibitem{label9} {P.G.~Grinevich, S.V.~Manakov}, {\it Inverse scattering problem for the two-dimensional Schr{\"o}dinger operator, the $ \bar{\partial} $-method and nonlinear equations}, Funct. Anal. Appl., 20 (1986), no. 2, 94–103.

\bibitem{label9_1} {R.G.~Novikov, G.M.~Henkin}, {\it The $ \bar{\partial} $-equation in the multidimensional inverse scattering problem}, Russian Math. Surveys, 42 (1987), no. 3, 109–180.

\bibitem{label10} {R.G.~Novikov}, {\it The inverse scattering problem on a fixed energy level for the two-dimensional Schr{\"o}dinger operator}, J. of Funct. Anal., 103 (1992), no. 2, 409–463.

\bibitem{label11} {R.G.~Novikov}, {\it Multidimensional inverse spectral problem for the equation $ - \Delta \psi + \left( v(x) - E u(x) \right) \psi = 0 $ }, Funct. Anal. Appl., 22 (1988), no. 4, 263-272.

\bibitem{label12} {R.G.~Novikov, M.~Santacesaria}, {\it Monochromatic reconstruction algorithms for twodimensional multi-channel inverse problems}, Int. Math. Res. Notices, 6 (2013), 1205-1229.

\bibitem{label13} {A.D.~Agaltsov, R.G.~Novikov}, {\it Riemann–Hilbert problem approach for two-dimensional flow inverse scattering}, J. Math. Phys., 55 (2014), no. 10, 103502.

\bibitem{label14} {V.A.~Burov, A.S.~Shurup, D.I.~Zotov, O.D.~Rumyantseva}, {\it Simulation of a functional solution to the acoustic tomography problem for data from quasi-point transducers}, Acoustical Physics, 59 (2013), no. 3, 345–360.

\bibitem{label15} {A.S.~Shurup, O.D.~Rumyantseva}, {\it Joint reconstruction of the speed of sound, absorption, and flows by the Novikov–Agaltsov functional algorithm}, Acoustical Physics, 63 (2017), no. 6, 751–768.

\bibitem{label16} {O.S.~Krasulin, A.S.~Shurup}, {\it Functional solution of ocean tomography problem with mode coupling}, Proc. of XVI Brekhovskikh’s Conference "Ocean Acoustics" and the XXXI session of the Russian Acoustical Society, M.: GEOS, 2018, 213-216.

\bibitem{label17} {M.I.~Belishev}, {\it Dynamical inverse problem for a Lam{\'e} type system}, J. of Inv. and Ill-posed Probl., 14 (2006), no. 8, 751-766.

\bibitem{label18} {M.I.~Belishev, A.L.~Pestov}, {\it Characterization of inverse data for one-dimensional two-velocity dynamical system} (in Russian), Algebra i Analiz, 26 (2014), no. 3, 89–130.

\bibitem{label19} {R.G.~Novikov}, {\it An iterative approach to non-overdetermined inverse scattering at fixed energy}, Sbornik: Mathematics, 206 (2015), no. 1, 120–134.

\bibitem{label20} {A.D.~Agaltsov, T.~Hohage, R.G.~Novikov}, {\it An iterative approach to monochromatic phaseless inverse scattering}, Inverse Problems, 35 (2019), no. 1, 024001.

\bibitem{label21} {R.G.~Novikov}, {\it Rapidly converging approximation in inverse quantum scattering in dimension 2}, Physics Letters A, 238 (1998), no. 2-3, 73–78.

\bibitem{label22} {R.G.~Novikov}, {\it Approximate inverse quantum scattering at fixed energy in dimension 2}, Proc. Steklov Inst. Math., 225 (1999), no. 2, 285–302.

\bibitem{label23} {V.A.~Burov, S.N. Vecherin, S.A. Morozov, O.D. Rumyantseva}, {\it Modeling of the exact solution of the inverse scattering problem by functional methods}, Acoustical Physics, 56 (2010), no. 4, 541–559.

\bibitem{label24} {V.A.~Burov, N.V. Alekseenko, O.D Rumyantseva}, {\it Multifrequency generalization of the Novikov algorithm for the two-dimensional inverse scattering problem}, Acoustical Physics, 55 (2009), no. 6, 843–856.

\bibitem{label25} {O.S.~Krasulin, A.S. Shurup}, {\it Numerical solution of three-dimensional problem of ocean adiabatic mode tomography based on functional-analytical algorithm}, Bull. Russ. Acad. Sci. Phys., 84 (2020), no. 2, 289-294.

\bibitem{label26} {S.I.~Kabanikhin, D.V.~Klyuchinskiy, N.S.~Novikov, M.A.~Shishlenin}, {\it Numerics of acoustical 2D tomography based on the conservation laws}, J. of Inv. and Ill-posed Probl., 28 (2020), no. 2, 287-297.

\bibitem{label27} {V.M.~Filatova, L.N.~Pestov, A.~Poddubskaya}, {\it Detection of velocity and attenuation inclusions in the medical ultrasound tomography}, J. of Inv. and Ill-posed Probl., 29 (2021), no. 3, 459-466.

\bibitem{label28} {A.V.~Goncharsky,  S.Y.~Romanov}, {\it Supercomputer technologies in inverse problems of ultrasound tomography}, Inverse Problems, 29 (2013), no. 7, 075004.

\bibitem{label29} {O.D.~Rumyantseva, A.S.~Shurup, D.I.~Zotov}, {\it Possibilities for separation of scalar and vector characteristics of acoustic scatterer in tomographic polychromatic regime}, J. of Inv. and Ill-posed Probl., 29 (2021), no. 3, 407-420.

\bibitem{label30} {O.D.~Rumyantseva, A.S.~Shurup}, {\it Equation for wave processes in inhomogeneous moving media and functional solution of the acoustic tomography problem based on it}, Acoustical Physics, 63 (2017), no. 1, 95–103.

\bibitem{label30_1} {A.I.~Nachman}, {\it Reconstruction from boundary measurements}, Annals of Math., 128 (1988), no. 3, 531–576.

\bibitem{label30_2} {Y.M.~Berezanskii}, {\it On the uniqueness theorem in the inverse problem of spectral analysis for the Schr{\"o}dinger equation}, Am. Math. Soc Trans., 35 (1964), 167–235.

\bibitem{label31} {V.A.~Burov,  M.N.~Rychagov}, {\it Diffraction tomography as inverse problem of scattering. interpolation approach. 1. Linearized version}, Acoustical Physics, 38 (1992), no. 4, 631-643.

\bibitem{label32} {A.A.~Goryunov, A.V. Saskovets}, {\it Inverse Scattering Problems in Acoustics} (in Russian), Moscow State University, Moscow, 1989.

\bibitem{label33} {V.A.~Burov, D.I.~Zotov, O.D.~Rumyantseva}, {\it Reconstruction of the sound velocity and absorption spatial distributions in soft biological tissue phantoms from experimental ultrasound tomography data}, Acoustical Physics, 61 (2015), no. 2, 231–248.


\end{thebibliography}
\end{document}